\documentstyle[12pt]{amsart}
\def\SBIMSMark#1#2#3{
 \font\SBF=cmss10 at 10 true pt
 \font\SBI=cmssi10 at 10 true pt
 \setbox0=\hbox{\SBF Stony Brook IMS Preprint \##1}
 \setbox2=\hbox to \wd0{\hfil \SBI #2}
 \setbox4=\hbox to \wd0{\hfil \SBI #3}
 \setbox6=\hbox to \wd0{\hss
             \vbox{\hsize=\wd0 \parskip=0pt \baselineskip=10 true pt
                   \copy0 \break%
                   \copy2 \break% 
                   \copy4 \break}}
 \dimen0=\ht6   \advance\dimen0 by \vsize \advance\dimen0 by 8 true pt
                \advance\dimen0 by -\pagetotal
 \dimen2=\hsize \advance\dimen2 by .25 true in
  \openin2=publishd.tex
  \ifeof2\setbox0=\hbox to 0pt{}
  \else 
     \setbox0=\hbox to 3.1 true in{
                \vbox to \ht6{\hsize=3 true in \parskip=0pt  \noindent  
                \input publishd.tex 
                \vfill}}
  \fi
  \closein2
  \ht0=0pt \dp0=0pt
 \ht6=0pt \dp6=0pt
 \setbox8=\vbox to \dimen0{\vfill \hbox to \dimen2{\copy0 \hss \copy6}}
 \ht8=0pt \dp8=0pt \wd8=0pt
 \copy8
 \message{*** Stony Brook IMS Preprint #1, #2 ***}
}

\renewcommand{\marginpar}[1]{}
\catcode`\@=12

\def\Empty{}
\newcommand\oplabel[1]{
  \def\OpArg{#1} \ifx \OpArg\Empty {} \else
  	\label{#1}
  \fi}
		
\long\def\realfig#1#2#3#4{
\begin{figure}[htp]
\centerline{\psfig{figure=#2,width=#4}}
\caption[#1]{#3}
\oplabel{#1}
\end{figure}}

\newcommand{\comm}[1]{}

\input{psfig}

\newcommand{\BBB}[1]{{\Bbb #1}}

\begin{document}
\SBIMSMark{1998/1a}{January 1998}{}
\title{On Biaccessible Points in the Julia set of a Cremer 
Quadratic Polynomial} 
\author{Dierk Schleicher and Saeed Zakeri} 
%\date{October 6, 1997}

\pagestyle{myheadings}
\markboth{\sc D. Schleicher and S. Zakeri}{\sc On Biaccessible Points}
\begin{abstract}
We prove that the only possible biaccessible points in the Julia set of a
Cremer quadratic polynomial are the Cremer fixed point and its preimages.  
This gives a partial answer to a question posed by C. McMullen on whether
such a Julia set can contain any biaccessible point at all. 
\end{abstract}
\maketitle

\noindent
{\bf \S 1. Introduction.} Every quadratic polynomial in the complex plane
with an irrationally indifferent fixed point is affinely conjugate to a
unique quadratic polynomial of the form  
\begin{equation}
P: z\mapsto e^{2\pi i \theta}z+z^2,
\end{equation}
where $0<\theta<1$ is irrational and the indifferent fixed point is located
at the origin. When $P$ is holomorphically linearizable about $0$, we call
it a {\it Siegel} polynomial. On the other hand, when the fixed point at the
origin is nonlinearizable, $P$ is called a {\it Cremer} polynomial. By the
theorem of Brjuno-Yoccoz \cite{Yoc}, $P$ is a Siegel polynomial if and only
if $\theta$ satisfies the {\it Brjuno condition}: 
$$\sum _{n=1}^{\infty} \frac{\log q_{n+1}}{q_n}< +\infty,$$
where the $q_n$ appear as the denominators of the rational approximations
coming from the continued fraction expansion of $\theta$. 

Recall that the {\it filled Julia set} of $P$ is
$$K(P)=\left\{ z\in {\BBB C}: \mbox{The orbit $\{ P^{\circ n}(z) \} _{n\geq
0}$ is bounded} \right\}$$ 
and the {\it Julia set} of $P$ is the topological boundary of the filled
Julia set: 
$$J(P)=\partial K(P).$$ 
Both sets are nonempty, compact, connected and the filled Julia set is full,
i.e., the complement $\BBB C \backslash K(P)$ is connected. Every connected
component of the interior of $K(P)$ is a topological disk called a {\it
bounded Fatou component} of $P$. According to Sullivan, every bounded Fatou
component must eventually map to the immediate basin of attraction of an
attracting periodic point, or to an attracting petal for a parabolic
periodic point, or to a periodic Siegel disk for $P$ (see for example
\cite{Mil}). On the other hand, by \cite{Douady1} a polynomial of degree
$d\geq 2$ can have at most $d-1$ nonrepelling periodic orbits. It follows
that in the Cremer case, $K(P)$ has no interior, so that $K(P)=J(P)$.\\ \\  
{\bf \S 2. Definitions.} Given a quadratic polynomial $P$, there exists a
conformal isomorphism 
$$\varphi: \overline{\BBB C} \backslash K(P)\rightarrow \overline{\BBB C}
\backslash \overline{\BBB D}$$ 
with $\varphi(\infty)=\infty$ and $\varphi'(\infty)>0$, which conjugates $P$
to the squaring map: 
\begin{equation}
\varphi(P(z))=(\varphi(z))^2.
\end{equation}
The $\varphi$-preimages of the radial lines and circles centered at the
origin are called the {\it external rays} and {\it equipotentials} of
$K(P)$. The external ray $R_t$, by definition, is  
$$\varphi^{-1} \{ re^{2 \pi i t}: r>1 \} ,$$
where $t\in \BBB R / \BBB Z$ is called the {\it angle} of the ray. From (2)
it follows that  
$$P(R_t)=R_{2t\ (mod\ 1)}.$$

We say that $R_t$ {\it lands} at $p\in J(P)$ if $\lim_{r\rightarrow 1}
\varphi^{-1} (re^{2 \pi i t})=p$. A point $p\in J(P)$ is called {\it
accessible} if there exists a simple arc in $\BBB C \backslash K(P)$ which
starts at infinity and terminates at $p$. According to a theorem of
Lindel\"{o}f (see for example \cite{Rudin}), $p$ is accessible exactly when
there exists an external ray landing at $p$. We call $p$ {\it biaccessible}
if it is accessible through at least two distinct external rays. By a
theorem of F. and M. Riesz \cite{Mil}, $K(P) \backslash \{ p \} $ is
disconnected whenever $p$ is biaccessible. It is interesting that the
converse is also true. More precisely, if there are at least $n>1$ connected
components of $K(P) \backslash \{ p \} $, then at least $n$ distinct
external rays land at $p$ (see for example \cite{McM}).\\ \\ 
{\bf \S 3. Some Known Results.} Very little is known about the topology of
the Julia set of $P$ in the Cremer case or the dynamics of $P$ on its Julia
set. The following theorem summarizes the basic results in this direction:\\
\\ 
{\bf Theorem 1.} {\it Let $P$ be a Cremer quadratic polynomial. Then} 
\begin{enumerate}
\item[(a)]
{\it The Julia set $J(P)$ cannot be locally connected} \cite{Sul}.
\item[(b)]
{\it Every neighborhood of the Cremer fixed point contains infinitely many 
repelling periodic orbits of $P$} \cite{Yoc}.
\item[(c)]
{\it The critical point $c=-e^{2 \pi i \theta}/2$ is recurrent, i.e., $c$
belongs to the closure of its orbit  $\{ P^{\circ n}(c) \} _{n> 0}$}
\cite{Man}. 
\item[(d)]
{\it The critical point $c$ is not accessible from $\BBB C \backslash J(P)$}
\cite{Kiw}.\\  
\end{enumerate}
See also \cite{Sor} for the so-called ``Douady's non-landing theorem'' which
partially explains why the Julia set of a generic Cremer quadratic
polynomial fails to be locally connected.\\ \\ 
{\bf \S 4. Hedgehogs.} Let $P$ be a Cremer quadratic polynomial. Let $U$ be
a simply connected domain with compact closure which contains the
indifferent fixed point $0$ of $P$. Suppose that $P$ is univalent in a
neighborhood of the closure $\overline U$. Then there exists a set $H=H_U$
with the following properties: 
\begin{enumerate}
\item[(i)]
$0\in H\subset \overline U \cap J(P),$
\item[(ii)]
$H$ is compact, connected and full,
\item[(iii)]
$H\cap \partial U$ is nonempty,
\item[(iv)]  
$P(H)=H.$
\end{enumerate}
Such an $H$ is called a {\it hedgehog} for the restriction
$P|_U:U\rightarrow \BBB C$. The existence of such completely invariant sets
is proved by Perez-Marco \cite{Per}.  

Hedgehogs turn out to be useful because of the following nice construction:
Uniformize the complement $\BBB C\backslash H$ by the Riemann map $\psi:\BBB
C \backslash H\rightarrow \BBB C \backslash \overline{\BBB D}$ and consider
the induced map $g=\psi\circ P \circ \psi^{-1}$ which is defined (by (iv)
above) and holomorphic in an open annulus $\{ z\in \BBB C: 1< |z| < r
\}$. Use the Schwarz Reflection Principle to extend $g$ to the annulus $\{
z\in \BBB C: r^{-1}< |z| < r \}$. The restriction of $g$ to the unit circle
$\BBB T$ will then be a real-analytic diffeomorphism whose rotation number
is exactly $\frac{1}{2 \pi i}\log P'(0)=\theta \in \BBB R/ \BBB Z$. This
allows us to transfer results from the more developed theory of circle
diffeomorphisms to the less explored theory of indifferent fixed points of
holomorphic maps.  

Using the above construction, it is not hard to prove the following fact 
(see \cite{Q} ):\\ \\
{\bf Proposition 2.} {\it The only point in a hedgehog which can be
biaccessible from outside of the hedgehog is the indifferent fixed point.}\\
\\ 
In fact, let $p$ be a biaccessible point in the hedgehog $H$ and $p\neq
0$. Then one can find a simple arc $\gamma$ in $\BBB C\backslash H$ which
starts and terminates at $p$ and does not encircle the indifferent fixed
point $0$. Let $D$ be the bounded connected component of $\BBB C\backslash
(H\cup \gamma)$ and let $D'=\psi(D)$. $D'$ is bounded by the simple arc
$\psi(\gamma)$ and an interval $I$ on the unit circle. Since $g$ has
irrational rotation number on the unit circle $\BBB T$, for some integer $N$
we have $\bigcup_{i=0}^N g^{\circ i}(I)=\BBB T$. By choosing $\gamma$ close
enough to $H$, we can assume that $g, g^{\circ 2}, \cdots, g^{\circ N}$ are
all defined on $D'$ and $\bigcup_{i=0}^N g^{\circ i}(D')$ contains an entire
outer neighborhood of $\BBB T$. It follows that $\bigcup_{i=0}^N P^{\circ
i}(D)$ covers an entire deleted neighborhood of $H$. This is possible only
if $0\in \overline D$ which contradicts our assumption. \\ \\ 
{\bf \S 5. The Theorem.} For the rest of this note, it will be more
convenient to use the normal form  
\begin{equation}
P: z\mapsto z^2+v
\end{equation}
for a quadratic polynomial $P$ in (1), which can be achieved by an affine
conjugation. In this normal form the origin is now the critical point,
$P(0)=v$ is the critical value, and the two fixed points are located at  
$$\alpha, \beta=\frac{1}{2}(1\pm \sqrt{1-4v}).$$
The indifferent fixed point $0$ in (1) is now located at the fixed point
$\alpha$, which is the fixed point further to the left (note that
$\alpha=\beta$ only when $v=1/4$).  

Since $P(z)=P(-z)$ by (3), the Julia set $J(P)$ is invariant under the
$180^{\circ}$ rotation $\tau: z\mapsto -z$. If $U$ is an open Jordan domain
in the plane such that $\overline{U}\cap \tau(\overline{U})=\emptyset$, it
follows that $P$ is univalent in some Jordan domain $V$ containing the
closure $\overline U$. 

Now we can state the main theorem:\\ \\
{\bf Theorem 3.} {\it Let $P:z\mapsto z^2+v$ be a Cremer quadratic
polynomial with the indifferent fixed point $\alpha$. Let $z\in J(P)$ be
biaccessible. Then there exists a non-negative integer $k$ such that
$P^{\circ k}(z)=\alpha$.}\\ \\ 
The proof begins as follows: Assume by way of contradiction that the forward
orbit of $z$ never hits the fixed point $\alpha$. Let $R$ and $R'$ be two
distinct external rays that land at $z$. Consider the connected component
$W$ of $\BBB C\backslash (R \cup R' \cup \{ z \} )$ which does not contain
$\alpha$. Since $\bigcup_{n\geq 0} P^{-n}(\alpha)$ is dense in $J(P)$ and
$W$ intersects the Julia set, there exists a {\it smallest} positive integer
$n$ such that $P^{-n}(\alpha)\cap W$ is nonempty. Let $\alpha_{-n}$ be a
point in this intersection. Define 
$$R_{n-1}=P^{\circ n-1}(R),\  R'_{n-1}=P^{\circ n-1}(R').$$ 
These rays are distinct because otherwise the forward orbit of $z$ would
have to hit the critical point $0$ which implies accessibility of $0$, which
contradicts Theorem 1(d). Clearly both rays land at $z_{n-1}=P^{\circ
n-1}(z)$. 

We claim that the ray pair $(R_{n-1}, R'_{n-1})$ separates $\alpha$ from its
preimage $-\alpha$. Assume it does not. Then the connected component
$W_{n-1}$ of $\BBB C\backslash (R_{n-1}\cup R'_{n-1}\cup \{ z_{n-1} \} )$
which contains  the fixed point $\alpha$ also contains $-\alpha$. Consider
$P^{\circ n-1}$  
as a finite-degree branched covering $\hat{\BBB C} \rightarrow \hat{\BBB C}$
and let $\{ O_1, \cdots , O_m \}$ be the collection of the connected
components of $P^{-(n-1)}(W_{n-1})$. Each $O_i$ is an unbounded simply
connected domain whose boundary is formed by two or more external rays in
$P^{-(n-1)}(R_{n-1}\cup R'_{n-1})$, and the mapping $P^{\circ n-1}:
O_i\rightarrow W_{n-1}$ is proper holomorphic. In particular, since
$W_{n-1}$ contains $\alpha$, it follows that each $O_i$ intersects
$P^{-(n-1)}(\alpha)$.  
On the other hand, $P^{\circ n-1}(\alpha_{-n})=-\alpha\in W_{n-1}$, which
means  
$\alpha_{-n}$ belongs to $O_i$ for some $1\leq i\leq m$. Since external rays
do not cross, it follows that $O_i \subset W$. But this means that $W$
intersects              $P^{-(n-1)}(\alpha)$, contradicting minimality of
$n$.  

Therefore, the ray pair $(R_{n-1}, R'_{n-1})$ separates the fixed point
$\alpha$ from its preimage $-\alpha$. Without loss of generality, we can
assume that these two rays actually separate $\alpha$ from the critical
point $0$, for otherwise we could replace them by the symmetric pair $(\tau
(R_{n-1}), \tau (R'_{n-1}))$. To simplify the notation, let $U=W_{n-1}$ be
the connected component of $\BBB C\backslash (R_{n-1}\cup R'_{n-1}\cup \{
z_{n-1} \} )$ which contains $\alpha$. To make $U$ have compact closure, one
can cut it off by some equipotential curve.  
%\\\\\\\\\\\\\\\\\\\\\\\\\\\\\
\realfig{qq}{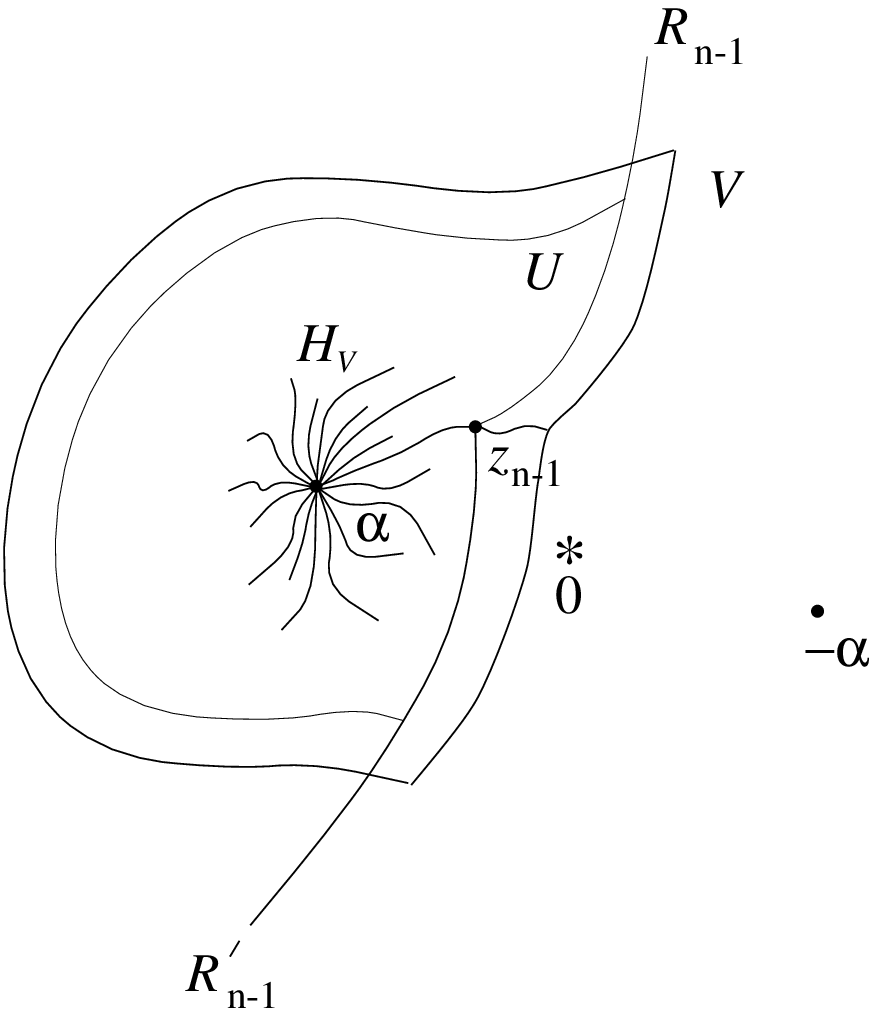}{}{7cm}
%\\\\\\\\\\\\\\\\\\\\\\\\\\\\\

First we observe that $P|_{\overline U}$ is univalent. In fact, if this were
not true, $\overline U$ would have to intersect $\tau (\overline U)$. Since
both $\overline U$ and $\tau (\overline U)$ are homeomorphic to the closed
disk, this would imply that the boundary ray pairs $(R_{n-1}, R'_{n-1})$ and
$(\tau (R_{n-1}), \tau (R'_{n-1}))$ intersect, which is impossible since
distinct rays are disjoint (in case $R'_{n-1}=\tau (R_{n-1})$, the landing
point would have to be the critical point which is ruled out by Theorem
1(d)). 

Now we consider a hedgehog $H_U$ for the restriction $P|_U: U\rightarrow
\BBB C$ as given in $\S 4$. Note that $z_{n-1}$ is the only point of the
Julia set on the boundary of $U$ and that $H_U\subset J(P)$ has to intersect
this boundary. Therefore, we simply have $z_{n-1}\in H_U$.  
By the remark before the statement of the theorem, we can find a slightly
larger Jordan domain $V \supset \overline U$ such that $P|_{\overline V}$ is
still univalent. The hedgehog $H_V$ for the restriction $P|_V: V\rightarrow
\BBB C$ has to contain $z_{n-1}$ also and reach the boundary of $V$. Since
$z_{n-1}$ is biaccessible from outside of the Julia set, it follows that
$H_V \backslash \{ z_{n-1} \}$ is disconnected. Therefore, $z_{n-1}$ is
biaccessible from outside of $H_V$. This contradicts Proposition 2, and
finishes the proof of Theorem 3.\\

\end{document}